\documentclass[12pt,twoside]{article}
\usepackage{graphicx,graphics, amsfonts, amsbsy,amssymb, amsmath, cite, array,arydshln, hyperref}
\usepackage{float}

\allowdisplaybreaks

\def\bpsp{\begin{pspicture}}
\def\epsp{\end{pspicture}}

\newtheorem{theorem}{Theorem}[section]

\newtheorem{lemma}[theorem]{Lemma}
\newtheorem{corollary}[theorem]{Corollary}

\newtheorem{conjecture}{Conjecture}

\usepackage{tikz}
\usetikzlibrary{chains,positioning}
\usepackage{amsmath,tikz,pgfplots}
\usetikzlibrary{decorations.pathmorphing}

\def\qed{\nolinebreak\hfill\rule{.2cm}{.2cm}\par\addvspace{.5cm}}

\begin{document}
\title{On distance Laplacian energy in terms of graph invariants}
\author{ Hilal A. Ganie$^a$, Rezwan Ul Shaban$^b$, Bilal A. Rather$^c$, S. Pirzada$^d$,  \\
$^{a}${\em Department of School Education, JK Govt. Kashmir, India }\\
$^{b,c,d}${\em Department of Mathematics, University of Kashmir, Srinagar, India}\\
$^a$hilahmad1119kt@gmail.com;~~$^b$rezwanbhat21@gmail.com;\\
 ~~$^c$bilalahmadrr@gmail.com; $^d$pirzadasd@kashmiruniversity.ac.in\\}

\date{}

\pagestyle{myheadings} \markboth{ Hilal, Rezwan, Bilal and Pirzada}{On distance Laplacian energy in terms of graph invariants}
\maketitle
\vskip 5mm
\noindent{\footnotesize \bf Abstract.} For a simple connected graph $G$ of order $n$ having distance Laplacian eigenvalues $ \rho^{L}_{1}\geq \rho^{L}_{2}\geq \cdots \geq \rho^{L}_{n}$, the distance Laplacian energy $DLE(G)$ is defined as $DLE(G)=\sum_{i=1}^{n}\left|\rho^{L}_i-\frac{2 W(G)}{n}\right|$, where $W(G)$ is the Wiener index of $G$.  We obtain a relationship between the Laplacian energy and distance Laplacian energy for graphs with diameter 2. We obtain lower bounds for the distance Laplacian energy $DLE(G)$ in terms of the order $n$, the Wiener index $W(G)$, independence number, vertex connectivity number and other given parameters. We characterize the extremal graphs attaining these bounds. We show that the complete bipartite  graph  has the minimum distance Laplacian energy among all connected bipartite graphs and complete split graph has the minimum distance Laplacian energy among all connected graphs with given independence number. Further, we obtain the distance Laplacian spectrum of the join  of a graph with the union of two other graphs. We show that the  graph $K_{k}\bigtriangledown(K_{t}\cup K_{n-k-t}), 1\leq t \leq \lfloor\frac{n-k}{2}\rfloor $,  has the minimum distance Laplacian energy among all connected graphs with vertex connectivity $k$. We conclude this paper with a discussion on trace norm of a matrix and the importance of our results in the theory of trace norm of the matrix $D^L(G)-\frac{2W(G)}{n}I_n$.
\vskip 3mm



\section{Introduction}
A graph is denoted by $G(V(G),E(G))$, where $V(G)=\{v_{1},v_{2},\ldots,v_{n}\}$ is its vertex set and $E(G)$ is its edge set. Throughout, $G$ is connected, simple and finite. $|V(G)|=n$ is the \textit{order} and $|E(G)|=m$ is the size of $G$. The set of vertices adjacent to $v\in V(G)$, denoted by $N(v)$, refers to the \textit{neighborhood} of $v.$ The \textit{degree} of $v,$ denoted by $d_{G}(v)$ (we simply write $d_v$ if it is clear from the context) means the cardinality of $N(v)$. A graph is \textit{regular} if each of its vertices have the same degree. The adjacency matrix $A=(a_{ij})$ of $G$ is a $(0, 1)$-square matrix of order $n$ whose $(i,j)$-entry is equal to $1$, if $v_i$ is adjacent to $v_j$ and equal to $0$, otherwise. Let $Deg(G)={diag}(d_1, d_2, \dots, d_n)$ be the diagonal matrix of vertex degrees $d_i=d_{G}(v_i)$, $i=1,2,\dots,n$ associated to $G$. The matrix $L(G)=Deg(G)-A(G)$ is the Laplacian matrix and its spectrum is the Laplacian spectrum of $G$. This matrix is real symmetric and positive semi-definite. We take
$0=\mu_n\leq\mu_{n-1}\leq\dots\leq\mu_1$ to be the Laplacian eigenvalues of $G$. The Laplacian energy of a graph \cite{IGBZ}, denoted by $LE(G)$, is defined as $LE(G) = \sum_{i=1}^n |\mu_i -\frac{2m}{n}|$. For some recent papers and related results on Laplacian energy, we refer to \cite{hbp,ph} and the references therein. For other undefined notations and terminology, the readers are referred to \cite{DMCM,lsg,sp}.\\
\indent The \textit{distance} between two vertices $u,v\in V(G),$ denoted by $d_{uv}$, is defined as the length of a shortest path between $u$ and $v$. The \textit{diameter} of $G$ is the maximum distance between any two vertices of $G.$ The \textit{distance matrix} of $G$, denoted by $D(G)$, is defined as $D(G)=(d_{uv})_{u,v\in V(G)}$. The \textit{transmission} $Tr_{G}(v)$ of a vertex $v$ is defined to be the sum of the distances from $v$ to all other vertices in $G$, i.e., $Tr_{G}(v)=\sum\limits_{u\in V(G)}d_{uv}.$ A graph $ G $ is said to be $k$-\textit{transmission regular} if $ Tr_{G}(v)=k,$ for each $ v\in V(G). $  The \textit{transmission} number or Wiener index of a graph $G$, denoted by $W(G), $ is the sum of distances between all unordered pairs of vertices in $G$.
Clearly, $W(G)=\frac{1}{2}\displaystyle\sum_{v\in V(G)}Tr_{G}(v)$. For any vertex $v_i\in V(G)$, the transmission $Tr_G(v_i)$ is called the \textit{transmission degree}, shortly denoted by $Tr_{i}$ and the sequence $\{Tr_{1},Tr_{2},\ldots,Tr_{n}\}$ is called the \textit{transmission degree sequence} of the graph $G$. \\
\indent Let $Tr(G)=diag (Tr_1,Tr_2,\ldots,Tr_n) $ be the diagonal matrix of vertex transmissions of $G$. Aouchiche and Hansen \cite{AH2} defined the \textit{distance Laplacian matrix} of $G$ as $ D^L(G)=Tr(G)-D(G) $.  \\
\indent Let $\rho^{D}_1\geq \rho^{D}_2\geq \dots\geq\rho^{D}_n$ and $\rho^{L}_1\geq \rho^{L}_2\geq \dots\geq\rho^{L}_n$ be respectively, the distance and distance Laplacian eigenvalues of the graph $G$. The distance energy \cite{IGV} of $G$ is the sum of the absolute values of the distance eigenvalues of $G$, that is, $DE(G)=\displaystyle\sum\limits_{i=1}^n|\rho^{D}_i|$. For some recent works on distance energy, we refer to \cite{dah,dr, AH1} and the references therein. \\
\indent The distance Laplacian energy $DLE(G)$ \cite{YYG} of a connected graph $G$ is defined as
\begin{equation*}
 DLE(G)=\sum\limits_{i=1}^n|\rho^{L}_i-\frac{2 W(G)}{n}|.
\end{equation*}
Let $\sigma$ be the largest positive integer such that $\rho^{L}_{\sigma}\geq \frac{2 W(G)}{n}$ and let $U_{k}^{L}(G)=\displaystyle\sum\limits_{i=1}^{k}\rho^{L}_i$ be the sum of $k$ largest distance Laplacian eigenvalues of $G$. Then using $\displaystyle\sum\limits_{i=1}^{n}\rho^{L}_i=2W(G)$, it can  be seen that
\begin{align*}
DLE(G)&=2\left(U_{\sigma}^L(G) -\frac{2\sigma W(G)}{n}\right) =2\max_{1\leq j\leq n}\left( \sum\limits_{i=1}^{j}\rho_{i}^L(G)-\frac{2jW(G)}{n} \right)\\&
=2\max_{1\leq j\leq n}\left(U_{j}^L(G)-\frac{2jW(G)}{n} \right).
\end{align*}
For some recent works on $DLE(G)$, see \cite{dah,dr,hilal}.\\
\indent The rest of the paper is organized as follows. In Section 2, we obtain a relationship between the Laplacian energy and distance Laplacian energy for graphs with diameter $2$. We also obtain a  lower bound for the distance Laplacian energy $DLE(G)$ in terms the order $n$, the Wiener index $W(G)$, etc. and characterize the extremal graphs attaining these bounds.  In  Section 3,  we study the distance Laplacian energy of connected bipartite graphs and of connected graphs with given independence number. We show that the complete bipartite  graph  has the minimum distance Laplacian energy among all connected bipartite graphs and complete split graph has the minimum distance Laplacian energy among all connected graphs with given independence number. In Section 4,  we study the distance Laplacian spectrum of join  of a graph with the union of two other graphs. We, also study the distance Laplacian energy of connected graphs with given vertex connectivity.  We show that the  graph $K_{k}\bigtriangledown(K_{t}\cup K_{n-k-t}), 1\leq t \leq \lfloor\frac{n-k}{2}\rfloor $,  has the minimum distance Laplacian energy among all connected graphs with vertex connectivity $k$. We have added a conclusion at the end to highlight the connection between the distance Laplacian energy of a connected graph $G$ with the trace norm of the matrix $D^L(G)-\frac{2W(G)}{n}I_n$.

\section{Bounds for the distance Laplacian energy of a graph}
We begin with the lemma, which gives the relation between the distance Laplacian spectrum of a graph and its connected spanning subgraph and  can be found in \cite{AH2}.
\begin{lemma}\label{lem1}
Let $ G $ be a connected graph of order $ n $ and size $ m $, where $ m\geq n $ and let $ G^{'}=G-e$ be a connected graph obtained from $ G $ by deleting an edge. Let  $ \rho_1^{L}(G)\geq  \rho_2^{L}(G)\geq \cdots\geq \rho_n^{L}(G) $ and $ \rho_1^{L}(G^{'})\geq \rho_2^{L}(G^{'})\geq \cdots\geq\rho_n^{L}(G^{'})  $ be respectively, the distance Laplacian eigenvalues of $ G $ and $ G^{'} $. Then $ \rho_i^{L}(G^{'})\geq \rho_i^{L}(G) $ holds for all $1\leq i\leq n . $
\end{lemma}

The following lemma shows that the distance Laplacian eigenvalues of a connected graph $G$ of diameter $2$ are completely determined by the Laplacian eigenvalues of the graph $G$ and can be found in \cite{AH2}.
\begin{lemma}\label{lem2}
Let $G$ be a connected graph of order $n\ge 2$ having diameter $d=2$. Let $\mu_1\ge \mu_2\ge \cdots\ge \mu_{n-1}> \mu_{n}=0$ be the Laplacian eigenvalues  and $\rho^{L}_1(G)\ge \rho^{L}_2(G)\ge\cdots\ge\rho^{L}_{n-1}(G)>\rho^{L}_n(G)=0$ be the distance Laplacian eigenvalues of $G$. Then $\rho^{L}_i(G)=2n-\mu_{n-i}$, for $i=1,2,\dots, n-1$.
\end{lemma}

A well known fact is that the complete graph $K_n$ is the only graph of diameter $1$. So for this graph, clearly $DLE(K_n)=2(n-1)$. Therefore, in the rest of the paper we will be dealing with graphs of diameter greater or equal to $2$.
The following result gives the distance Laplacian energy of a graph of diameter $2$ in terms of sum of the Laplacian eigenvalues of $G$.
\begin{theorem}\label{thm21}
Let $G$ be a connected graph of order $n\ge 3$ and size $m$ having diameter  $2$. Then
\begin{align*}
DLE(G)=2\Big(\sigma\Big(\frac{2m}{n}+2\Big)-2m+S_{n-\sigma-1}(G)\Big),
\end{align*} where $S_{n-\sigma-1}(G)=\sum\limits_{i=1}^{n-\sigma-1}\mu_{i}$ is the sum of the $n-\sigma-1$ largest Laplacian eigenvalues of $G$ and $\sigma, 1\leq \sigma\le n-2$, is  the number of distance Laplacian eigenvalues of $G$ which are greater than or equal to $\frac{2W(G)}{n}$.
\end{theorem}
{\bf Proof.}
Let $G$ be a connected graph of order $n$ having $m$ edges. Since diameter of $G$ is two, it follows that $Tr(v_i)=d_i+2(n-1-d_i)=2n-2-d_i$, for all $v_i\in V(G)$ and so $2W(G)=2n(n-1)-2m$. Let $\sigma$ be the number of distance Laplacian eigenvalues of $G$ which are greater than or equal to $\frac{2W(G)}{n}$. Using Lemma \ref{lem2} and  the definition of Laplacian energy, we have
\begin{align*}
DLE(G)&=2\left( \sum\limits_{i=1}^{\sigma}\rho_i^L(G)-\frac{2\sigma W(G)}{n} \right)\\&
= 2\left( \sum\limits_{i=1}^{\sigma}(2n-\mu_{n-i})-\sigma(2n-2-\frac{2m}{n}) \right)\\&
= 2\left( \sigma \big(\frac{2m}{n}+2\big)-\sum\limits_{i=1}^{\sigma}\mu_{n-i} \right)
=2\Big(\sigma\big(\frac{2m}{n}+2\big)-2m+S_{n-\sigma-1}(G)\Big),
\end{align*}
where $S_{n-\sigma-1}(G)=\sum\limits_{i=1}^{n-\sigma-1}\mu_{i}$ is the sum of the $n-\sigma-1$ largest Laplacian eigenvalues of $G$.\qed

From Theorem \ref{thm21}, it is clear that any lower bound or upper bound for the graph invariant $S_k(G)$, the sum of $ k $ largest Laplacian eigenvalues of $ G $, gives a lower bound or  upper bound for $DLE(G)$. In fact, there is a conjecture by A. Brouwer for  the  graph invariant $S_k(G)$, which is stated as follows.
\begin{conjecture}\label{Conjecture 1.2.} If $G$ is any graph with order $n$ and size
$m$, then
\[S_k(G)=\sum\limits_{i=1}^{k}\mu_i\leq m+\binom{k+1}{2}, \qquad \text{for any} \qquad k \in \{1,2,\dots,n\}.\]
\end{conjecture}
This conjecture has been shown to be true for various families of graphs, but as a whole this conjecture is still open. For some recent developments on Brouwer's conjecture, we refer to \cite{hpbt} and the references therein.

The following theorem gives a relation between the distance Laplacian energy and Laplacian energy of a graph of diameter $2$.
\begin{theorem}\label{thm23}
Let $G$ be a connected graph of order $n\ge 3$ and size $m$ having diameter $2$. Let $\sigma$ and $t$, $1\le \sigma, t\le n-2$, be respectively the number of distance Laplacian eigenvalues and number of Laplacian eigenvalues
of $G$ which are greater than or equal to $\frac{2W(G)}{n}$ and $\frac{2m}{n}$. Then
\begin{align*}
LE(G)-2\Big(\frac{2m}{n}-2(n-1)+2t\Big)\leq DLE(G)\le LE(G)+4\Big(\sigma-\frac{m}{n}\Big).
\end{align*}
\end{theorem}
 {\bf Proof.}
 Let $\sigma$ be the number of distance Laplacian eigenvalues of $G$ which are greater than or equal to $\frac{2W(G)}{n}$. Then by definition of distance Laplacian energy, we have
 \begin{align}\label{a}
 DLE(G)=2\Big(U^{L}_{\sigma}(G)-\frac{2\sigma W(G)}{n}\Big)=2\max_{1\le j\le n-1 }\Big(U^{L}_{j}(G)-\frac{2j W(G)}{n}\Big).
 \end{align}
 Also, if $t$ is the number of Laplacian eigenvalues of $G$ which are greater than or equal to $\frac{2m}{n}$, then by definition of Laplacian energy, we have
 \begin{align}\label{b}
  LE(G)=2\Big(S_{t}(G)-\frac{2t m}{n}\Big)=2\max_{1\le j\le n-1 }\Big(S_{j}(G)-\frac{2j m}{n}\Big).
  \end{align}
 Using Theorem \ref{thm21} and second equality of (\ref{b}), we have
 \begin{align*}
 DLE(G) =&2\Big(\sigma\big(\frac{2m}{n}+2\big)-2m+S_{n-\sigma-1}(G)\Big)\\&
= 2\Big(2\sigma-\frac{2m}{n}\Big)+2\Big(S_{n-\sigma-1}(G)-\frac{2m(n-\sigma-1)}{n}\Big)\\&
 \le LE(G)+4\Big(\sigma-\frac{m}{n}\Big),
 \end{align*} as $1\le n-\sigma-1\le n-1$.
 Using Lemma \ref{lem2} and (\ref{a}), we have
 \begin{align*}
   LE(G)=&2\Big(S_{t}(G)-\frac{2t m}{n}\Big)=2\Big((2n-\frac{2m}{n})t-\sum\limits_{i=1}^{t}\rho^{L}_{n-i}\Big)\\&
   =2\Big((2n-\frac{2m}{n})t-2W(G)+\sum\limits_{i=1}^{n-t-1}\rho^{L}_{i}\Big)\\&
   =2\Big(\frac{2m}{n}-2(n-1)+2t\Big)+2\Big(\sum\limits_{i=1}^{n-t-1}\rho^{L}_{i}-\frac{2(n-t-1)W(G)}{n}\Big)\\&
   \le 2\Big(\frac{2m}{n}-2(n-1)+2t\Big)+DLE(G),
   \end{align*}
 as $2W(G)=2n(n-1)-2m$ and $1\leq n-t-1\le n-1$. This completes the proof.\qed

From Theorem \ref{thm23}, the following observation is immediate.
\begin{corollary}
Let $G$ be a connected graph of order $n\geq 3$ and size $m$ having diameter $2$. Let $\sigma$ and $t$, $1\le \sigma, t\le n-2$, be respectively the number of distance Laplacian eigenvalues and number of Laplacian eigenvalues
of $G$ which are greater than or equal to $\frac{2W(G)}{n}$ and $\frac{2m}{n}$.  Then $LE(G)>DLE(G)$, provided that $m>\sigma n$ and $LE(G)< DLE(G)$, provided that $t<n-\frac{m}{n}-2$.
\end{corollary}
 For the star graph $K_{1,n-1}$  the Laplacian spectrum is $\{n,1^{[n-2]},0\}$ and $m=n-1$. It is easy to see that $t=1$ for $K_{1,n-1}$ and so the inequality $LE(K_{1,n-1})< DLE(K_{1,n-1})$ is valid for $ n\geq 4 $.

 \indent The following observation also follows from Theorem \ref{thm23}.
 \begin{corollary}
Let $G$ be a connected graph of order $n\geq 3$ and size $m$ having diameter $2$. Let $\sigma$ and $t$, $1\le \sigma, t\le n-2$, be respectively the number of distance Laplacian eigenvalues and number of Laplacian eigenvalues
of $G$ which are greater than $\frac{2W(G)}{n}$ and $\frac{2m}{n}$. Then $\sigma\geq n-(t+1)$.
 \end{corollary}

 The following lemma \cite{dah} gives an upper bound for the second smallest distance Laplacian eigenvalue $\rho^{L}_{n-1}(G)$ in terms of order $n$ and the minimum transmission degree $Tr_{\min}$ of the graph $G$.
\begin{lemma}\label{lem3}
Let $G$ be a connected graph of order $n\geq 3$ having minimum transmission $Tr_{\min}$ and second smallest distance Laplacian eigenvalue $\rho^{L}_{n-1}(G)$. Then
\begin{align}\label{1}
\rho^{L}_{n-1}(G)\leq \frac{n}{n-1}Tr_{\min},
\end{align}
with equality if and only if $G$ is a graph containing a vertex of transmission degree $n-1$.
\end{lemma}

Now, we obtain a lower bound for the distance Laplacian energy in terms the order $n$, the Wiener index $W(G)$ and the minimum transmission $Tr_{\min}$ of the graph $G$.
\begin{theorem}\label{thm2}
Let $ G $ be a connected graph of order $ n\geq 3 $ having minimum transmission degree $Tr_{\min}$ and Wiener index $W(G)$. Then
\begin{align}\label{4}
DLE(G) \geq \frac{8W(G)}{n}-\frac{2n Tr_{\min}}{n-1},
\end{align}
equality occurs if and only if $\sigma(G)=n-2$ and $ G$  containing a vertex of transmission degree $n-1$.
\end{theorem}
{\bf Proof.} Let $G$ be a connected graph of order $n\geq 3$ having distance Laplacian eigenvalues $\rho^{L}_{1}(G)\geq \rho^{L}_{2}(G)\geq\cdots\geq\rho^{L}_{n-1}(G)\geq \rho^{L}_{n}(G)=0$. Let $\sigma$ be the number of distance Laplacian eigenvalues of $G$ which are greater than or equal to $\frac{2W(G)}{n}$. Using $\sum\limits_{i=1}^{n}\rho^{L}_i(G)=2W(G)$ and the definition of Laplacian energy \cite{dah}, we have
\begin{align*}
DLE(G)&=2\max_{1\leq j\leq n-1}\left( \sum\limits_{i=1}^{j}\rho_i^L(G)-\frac{2jW(G)}{n} \right)\\&
\geq 2\left( \sum\limits_{i=1}^{n-2}\rho^{L}_i(G)-\frac{2(n-2)W(G)}{n} \right)\\&
=2\left(2W(G)-\rho^{L}_{n-1}(G)-\frac{2(n-2)W(G)}{n}\right)=\frac{8W(G)}{n}-2\rho^{L}_{n-1}(G)\\&
\geq \frac{8W(G)}{n}-\frac{2n Tr_{\min}}{n-1}.
\end{align*}
Equality occurs in (\ref{4}) if and only if equality occurs in
\begin{align}\label{5}
\max_{1\leq j\leq n-1}\left( \sum\limits_{i=1}^{j}\rho_i^L(G)-\frac{2jW(G)}{n} \right)=\left( \sum\limits_{i=1}^{n-2}\rho^{L}_i(G)-\frac{2(n-2)W(G)}{n} \right)
\end{align}
 and equality occurs in $\rho^{L}_{n-1}(G)\leq \frac{n}{n-1}Tr_{\min}$. Equality occurs in (\ref{5}) if and only if $\sigma=n-2$ and by Lemma \ref{lem3} equality occurs in $\rho^{L}_{n-1}(G)\leq \frac{n}{n-1}Tr_{\min}$ if and only if $G$ contains a vertex having transmission degree $n-1$. This shows that Equality occurs in (\ref{4}) if and only if  $\sigma=n-2$ and $G$ contains a vertex having transmission degree $n-1$. This completes the proof.\qed

The parameter $t$ gives the number of Laplacian eigenvalues of  a graph $G$ which are in the interval $[\frac{2m}{n},n]$. This parameter has been studied for various families  of graphs and is presently as active topic of research in the field of spectral theory of graphs. Like the parameter $t$,
the parameter $\sigma$ is concerned with the distribution of distance Laplacian eigenvalues of a connected graph $G$. In fact, the value of $\sigma$ gives the number of distance Laplacian eigenvalues which are in the interval $[\frac{2W(G)}{n},\infty)$. It will be of interest to discuss the following problem for the parameter $\sigma$.\\
{\bf Problem 1.} Characterize the graphs having $\sigma=1$, $\sigma=n-2$ and $\sigma=n-1$. Establish relations between $\sigma$ with different parameters of the graph $G$.

For the transmission regular graphs $G$, it can be seen that $\sigma= n-\gamma$, where $\gamma$ is the positive inertia of the distance matrix of the  graph $G$.

\section{Distance Laplacian energy of  bipartite graphs and graphs with given independence number}

In this section, we show that among all connected bipartite graphs the complete bipartite graph has the minimum distance Laplacian energy. We also show that among all connected graphs with independence number $t$, the complete split graph has the minimum distance Laplacian energy.

Among all connected bipartite graphs, the complete bipartite graph has the minimum distance Laplacian energy can be seen as follows.
\begin{theorem}\label{thm31}
Let $ G $ be a connected bipartite graph of order $ n\geq 3 $ with partite sets of cardinality $ a $ and $ b $, $a+b=n$.\\
{\bf (1).} If $ a<b $ or $ a>b $, then
\begin{align*}
DLE(G)\geq\left\{\begin{array}{lr}4n^2-6n-4ab-\frac{4(n-2)W(G)}{n}, &\mbox{if $2ab\geq n(b-2)\ or\ 2ab\geq n(a-2)$},\\
2(b-1)\left( 2n-a-\dfrac{2W(G)}{n} \right), &\mbox{if $2ab < n(b-2),~a<b$},\\
2(a-1)\left( 2n-b-\dfrac{2W(G)}{n} \right), &\mbox{if $2ab < n(a-2),~a>b$},
\end{array} \right.
\end{align*}
with equality if and only if $ G\cong K_{a,b} $.\\
{\bf (2).} If $ n=2a $ and $ n\geq 5 ,$ then $ DLE(G)\geq 12a(a-1)-\frac{4(n-2)W(G)}{n}=12b(b-1)-\frac{4(n-2)W(G)}{n}$. Equality occurs  if and only if $ G\cong K_{a,a}.$
\end{theorem}
\noindent\textbf{Proof.} Let $ G $ be a connected bipartite graph of order $ n $ with vertex set $ V(G) .$ Let $ V(G)=V_1\cup V_2 ,$ with $ |V_1|=a, |V_2|=b $ be a bipartition of the vertex set $ V(G)$ of $G$. Since  $ G $ is a connected bipartite  graph with partite sets of cardinality $a$ and $b$, it follows that $G$ is a spanning subgraph of the complete bipartite graph $ K_{a,b}$. Therefore, by Lemma \ref{lem1}, we have $ \partial_i^L(G)\geq\partial_i^L(K_{a,b})$, for all $i=1,2,\dots,n$. With this and the definition of distance Laplacian energy $DLE(G)$ , we have
\begin{align}\label{31}
DLE(G)&=2\left(\sum\limits_{i=1}^{\sigma}\rho_i^L(G) -\frac{2\sigma W(G)}{n}\right) =2\max_{1\leq j\leq n}\left( \sum\limits_{i=1}^{j}\rho_i^L(G)-\frac{2jW(G)}{n} \right)\nonumber\\&
\geq 2\max_{1\leq j\leq n}\left( \sum\limits_{i=1}^{j}\rho_i^L(K_{a,b})-\frac{2jW(G)}{n} \right),
\end{align}
 where $ \sigma $ is the largest positive integer such that $ \rho_{\sigma}^L(G)\geq \frac{2W(G)}{n} .$\\
\indent In \cite{AH}, it can be seen that the distance Laplacian spectrum  of $ K_{a,b}$ is $ \{ (2n-a)^{[b-1]}, (2n-b)^{[a-1]},n,0\}$ and  $2W(K_{a,b})=2n^2-2n-2ab $. Since $a$ and $b$ are positive integers, therefore by trichotomy law, we have $a<b$ or $a>b$ or $a=b$. We first prove the result for $ a<b $. If $ a<b $, then $ 2n-a\geq 2n-b $. Also $ 2n-a>n$ holds as $ n>a $. Likewise, if $ n>b $, then $ 2n-b>n $. These observations imply that $ 2n-a $ is the distance Laplacian spectral radius of $ K_{a,b} $ in case of $ a<b $ and so we always have $ 2n-a\geq \frac{2W(K_{a,b})}{n}$. The case $ a>b $ can be discussed similarly. For the eigenvalue $n$, we have $n<\frac{2W(K_{a,b})}{n}$ implying that $n^2-2n-2ab>0$, which further gives
\begin{align}\label{32}
2a^2-2an+n^2-2n>0, \quad \text{as}\quad  a+b=n.
\end{align} Consider the polynomial $f(a)=2a^2-2an+n^2-2n$, for $1\leq a<n$. The discriminant of this polynomial is $d=4n(4-n)$. Clearly for $n>4$, we have $d<0$ and so for this $n$, we always have $f(a)>0$. For $n=3,4$, it can be seen by direct calculation  that inequality (\ref{32}) holds. This shows that $n<\frac{2W(K_{a,b})}{n}$ holds for all $a<b$.
 For the eigenvalue $ 2n-b ,$ we have $2n-b\geq \frac{2W(K_{a,b})}{n}$  implying that $2ab\geq n(b-2)$. This shows that if  $ 2ab\geq n(b-2) $, then $ \frac{2W(G)}{n}\leq 2n-b $ and if $ 2ab< n(b-2) $, we have $ \frac{2W(G)}{n}> 2n-b.$ In other words, if $ 2ab\geq n(b-2)$, then $\sigma=n-2$ while if $ 2ab< n(b-2) $, then $\sigma=b-1$.  Therefore, we have the following cases to consider.\\
{\bf Case (i).} If $ 2ab\geq n(b-2)$, then the  number of eigenvalues which are greater or equal to $ \frac{2W(G)}{n} $ are $ n-2 $, that is $ \sigma=n-2 .$ $ Since  $$1\leq \sigma\leq n-1$, from \eqref{31} it follows that
\begin{align*}
DLE(G) & \geq 2\left( \sum\limits_{i=1}^{n-2}\rho_i^L(K_{a,b})-\frac{2(n-2)W(G)}{n}\right)\\
& = 2\left((b-1)(2n-a)+(a-1)(2n-b)- \frac{2(n-2)W(G)}{n}\right)\\
&= 4n^2-6n-4ab-\frac{4(n-2)W(G)}{n}.
\end{align*}
{\bf Case(ii).} If $ 2ab< n(b-2) $,  then the  number of eigenvalues, which are greater or equal to $ \frac{2W(G)}{n} $, is $ b-1$,that is $ \sigma=b-1.$ Since $1\leq \sigma\leq n-1$, from \eqref{31} it follows that
\begin{align*}
DLE(G) & \geq 2\left(\sum\limits_{i=1}^{b-1}\rho_i^L(K_{a,b})-\frac{2(b-1)W(G)}{n}\right)\\
& = 2(b-1)\left( 2n-a-\dfrac{2W(G)}{n} \right).
\end{align*}
If $a>b$, then proceeding similarly as above, it can be seen that $2n-b$ is the distance Laplacian spectral radius, $n<\frac{2W(K_{a,b})}{n}$ holds for all $a,b$, while as $ \frac{2W(G)}{n}\leq 2n-a $, for $ 2ab\geq n(a-2) $ and $ \frac{2W(G)}{n}> 2n-a,$ for $ 2ab< n(a-2)$. Using this information in  \eqref{31}, we get
\begin{align*}
DLE(G)\geq 4n^2-6n-4ab-\frac{4(n-2)W(G)}{n},
\end{align*}
for $ 2ab\geq n(a-2)$ and
\begin{align*}
DLE(G)\geq 2(a-1)\left( 2n-b-\dfrac{2W(G)}{n} \right),
\end{align*}
for $ 2ab< n(a-2)$. Equality occurs in each of the inequalities above if and only if equality occurs in \eqref{31}. It is clear that equality occurs in \eqref{31} if and only if $G\cong K_{a,b}$. This gives that equality occurs if and only if $G\cong K_{a,b}$. This completes the proof for these cases.\\
\indent If $a=b$, then  $ n=2a $ and so the distance Laplacian spectrum of $ K_{a,a}$ is $\{ (2n-a)^{[n-2]},n,0 \}$ and $ \dfrac{2W(K_{a,a})}{n}=\dfrac{2n^2-an-3n+2a}{n}.$ Since $ 2n-a \geq \dfrac{2W(K_{a,a})}{n}$ holds for all $a$ and for $ n\geq 5$, we have $ \dfrac{2W(K_{a,a})}{n}>n$, it follows that $ \sigma =n-2$. Therefore, from \eqref{31}, we have
\begin{align*}
DLE(G)&\geq 2\left(\sum\limits_{i=1}^{n-2}\rho_i^L(K_{a,b})-\frac{2(n-2)W(G)}{2} \right)\\&
 = 12a(a-1)-\frac{4(n-2)W(G)}{n}.
\end{align*}
Equality case can be discussed similarly as above. This completes the proof.\qed

For the complete bipartite graph $K_{a,b}$ with $a+b=n$ and $n\geq 3$, using Theorem \ref{thm31}, we have the following:\\
\indent  For $2ab\geq n(b-2)$ with $a<b$; or $2ab\geq n(a-2)$ with $ a>b$, we have
$4n^2-6n-4ab-\frac{4(n-2)W(K_{a,b})}{n}>  12a(a-1)-\frac{4(n-2)W(K_{a,b})}{n}$ implying that $8a^2+(4n-12)a-(4n^2-6n)<0$. The zeros of the polynomial $f(a)=8a^2+(4n-12)a-(4n^2-6n)$ are $-n+\frac{3}{2}, \frac{n}{2}$. This implies that $f(a)<0$, for all $a\in (-n+\frac{3}{2},\frac{n}{2})$. This shows that $4n^2-6n-4ab-\frac{4(n-2)W(K_{a,b})}{n}>  12a(a-1)-\frac{4(n-2)W(G)}{n}$, holds for all $a<\frac{n}{2}$. For $2ab < n(b-2),~a<b$, we have $2(b-1)\left( 2n-a-\dfrac{2W(K_{a,b})}{n} \right)>12a(b-1)-\frac{4(2b-2)W(K_{a,b})}{n}$ giving that $2n-a+\dfrac{2W(K_{a,b})}{n}>6a$. Using $2W(K_{a,b})=2n^2-2n-2an+2a^2=(n-1)^2+(n-a)^2+a^2-1$, we get $2a^2-9an+4n^2-2n>0$. The zeros of the polynomial $g(a)=2a^2-9an+4n^2-2n$ are $y_1=\frac{9n+\sqrt{49n^2+16n}}{4}, y_2=\frac{9n-\sqrt{49n^2+16n}}{4}$. This shows that $g(a)>0$, for all $a>y_1$ and for all $a<y_2$. Since $y_1>n$, it follows that $g(a)>0$, for all $a<y_2$. For $n\geq 3$, it is easy to see that $y_2>\frac{n-1}{2}$, implying that $g(a)>0$, for all $a<\frac{n-1}{2}$. This shows that $2(b-1)\left( 2n-a-\dfrac{2W(K_{a,b})}{n} \right)>12a(b-1)-\frac{4(2b-2)W(K_{a,b})}{n}$, holds for all  $a<\frac{n-1}{2}$. For $2ab < n(a-2),~a>b$, proceeding similarly as above it can be seen that $2(a-1)\left( 2n-b-\dfrac{2W(K_{a,b})}{n} \right)>12b(a-1)-\frac{4(2a-2)W(K_{a,b})}{n}$, holds for all  $b<\frac{n-1}{2}$. Now, using Theorem \ref{thm31} and the fact that $a$ and $b=n-a$ are positive integers, we have the following observation.
\begin{corollary}
Among all bipartite graphs of order $n\geq 3$, the  complete bipartite graph $K_{\lfloor\frac{n}{2}\rfloor,\lceil\frac{n}{2}\rceil}$ has the minimum distance Laplacian energy.
\end{corollary}

The next observation follows from Lemma \ref{lem1} and Theorem \ref{thm31}.
\begin{corollary}
Let $G$ be a   connected bipartite graph with partite sets of cardinality $a$ and $b$, $a+b=n$. Let $\rho_{1}^{L}(G)\ge \rho_{2}^{L}(G)\cdots\ge \rho_{n-1}^{L}(G)>\rho_{n}^{L}(G)=0$ be the distance Laplacian eigenvalues of $G$. Then
$\rho_{i}^{L}(G)\geq 2n-a,$ for all $1\leq i\leq b-1$, $\rho_{i}^{L}(G)\geq 2n-b,$ for all $b\leq i\leq a+b-2$, $\rho_{n-1}^{L}(G)\geq n$. Equality occurs in each of these inequalities if and only if $G\cong K_{a,b}$.
\end{corollary}

A complete split graph, denoted by $CS_{t,n-t}$, is the graph consisting of a clique on $t$ vertices and an independent set (a subset of vertices of a graph is said to be an independent set if the subgraph induced by them is an empty graph) on the remaining $n-t$ vertices, such that each vertex of the clique is adjacent to every vertex of the independent set.

The following theorem shows that among all connected graphs with given independence number $n-t$, $1\leq t\leq n-1$, the complete split graph $CS_{t,n-t}$ has the minimum distance Laplacian energy.
\begin{theorem}\label{thm3}
Let $ G $ be a connected graph of order $n\geq 3$ having independence number $n-t$, $1\leq t\leq n-1$. Then
\begin{align*}
DLE(G)\geq\left\{\begin{array}{lr}2(n-t-1)\Big(2n-t-\frac{2W(G)}{n}\Big), &\mbox{if $t<n-\frac{1}{2}-\sqrt{n+\frac{1}{4}}$},\\
4n^2-2nt-4n+2t(t+1)-\frac{4(n-1)W(G)}{n}, &\mbox{if $t\geq n-\frac{1}{2}-\sqrt{n+\frac{1}{4}}$},
\end{array} \right.
\end{align*}
equality occurs in each case if and only if $ G\cong CS_{t,n-t}, ~1\leq t\leq n-1$.\\
\end{theorem}
\noindent\textbf{Proof.} Let $ G $ be a connected graph of order $ n\geq  3$ having independence number $n-t$. Let $CS_{t,n-t}$ be the complete split graph having independence number $n-t$. It is clear that $G$ is a spanning subgraph of $CS_{t,n-t}$. Therefore, by Lemma \ref{lem1}, we have $ \rho_i^L(G)\geq\rho_i^L(CS_{t,n-t})$. Let $ \sigma $ be the largest positive integer such that $ \rho_{\sigma}^L(G)\geq \frac{2W(G)}{n}.$ With this information, it follows that
\begin{align}\label{eq2}
DLE(G)&=2\left(\sum\limits_{i=1}^{\sigma}\rho_i^L(G) -\frac{2\sigma W(G)}{n}\right) =2\max_{1\leq j\leq n}\left( \sum\limits_{i=1}^{j}\rho_i^L(G)-\frac{2jW(G)}{n} \right)\nonumber\\&
\geq 2\max_{1\leq j\leq n}\left( \sum\limits_{i=1}^{j}\rho_i^L(CS_{t,n-t})-\frac{2jW(G)}{n} \right).
\end{align}
The distance Laplacian spectrum \cite{AH} of the complete split graph $ CS_{t,n-t} $ is $\{(2n-t)^{[n-t-1]},n^{[t]},0\}$ with $ \frac{2W(CS_{t,n-t})}{n}=\frac{2n(n-t-1)+t(t+1)}{n}$. Since $n-t\geq 1$, it follows that $2n-t$ is the distance Laplacian spectral radius of $CS_{t,n-t}$. For the eigenvalue $n$, we have $n< \frac{2W(CS_{t,n-t})}{n}= \frac{2n(n-t-1)+t(t+1)}{n}$, which after simplification gives
 \begin{align}\label{7}
 t^2-(2n-1)t+n^2-2n>0.
 \end{align}
 Consider the polynomial $f(t)=t^2-(2n-1)t+n^2-2n$, for $1\leq t\leq n-1$. The roots of this polynomial are $x_1=n-\frac{1}{2}-\sqrt{n+\frac{1}{4}}$ and $x_2= n-\frac{1}{2}+\sqrt{n+\frac{1}{4}}$. This implies that $f(t)>0$ for all $t<x_1$ and $f(t)>0$ for all $t>x_2$. Since $x_2>n$ and $t\leq n-1$, it follows that inequality (\ref{7}) holds for all $t<n-\frac{1}{2}-\sqrt{n+\frac{1}{4}}$. From this, it follows that for $t<n-\frac{1}{2}-\sqrt{n+\frac{1}{4}}$, we have $\sigma=n-t-1$ and for $t\geq n-\frac{1}{2}-\sqrt{n+\frac{1}{4}}$, we have $\sigma=n-1$. Therefore, we have the following cases to consider.\\
{\bf Case (i).} If $t<n-\frac{1}{2}-\sqrt{n+\frac{1}{4}}$, then $\sigma=n-t-1$. Since $1\leq \sigma\leq n-1$, from (\ref{eq2}) it follows that
\begin{align*}
DLE(G) &\geq 2\max_{1\leq j\leq n}\left( \sum\limits_{i=1}^{j}\rho_i^L(CS_{t,n-t})-\frac{2jW(G)}{n} \right)\\&
\geq 2\left( \sum\limits_{i=1}^{n-t-1}\rho_i^L(CS_{t,n-t})-\frac{2(n-t-1)W(G)}{n} \right)\\&
=2(n-t-1)\Big(2n-t-\frac{2W(G)}{n}\Big).
\end{align*}
{\bf Case (ii).} If $t\geq n-\frac{1}{2}-\sqrt{n+\frac{1}{4}}$, then $\sigma=n-1$. Since $1\leq \sigma\leq n-1$, from (\ref{eq2}) it follows that
\begin{align*}
DLE(G) &\geq 2\max_{1\leq j\leq n}\left( \sum\limits_{i=1}^{j}\rho_i^L(CS_{t,n-t})-\frac{2jW(G)}{n} \right)\\&
\geq 2\left( \sum\limits_{i=1}^{n-1}\rho_i^L(CS_{t,n-t})-\frac{2(n-1)W(G)}{n} \right)\\&
=2(2n-t)(n-t-1)+2nt-\frac{4(n-1)W(G)}{n}\\&
=4n^2-2nt-4n+2t(t+1)-\frac{4(n-1)W(G)}{n}.
\end{align*}
Equality occurs in all the inequalities above if and only if equality occurs in (\ref{eq2}). It is clear that equality occurs in (\ref{eq2}) if and only if $G\cong  CS_{t,n-t}$. This shows that equality occurs in all the inequalities above if and only if $G\cong  CS_{t,n-t}$. This completes the proof.\qed

We characterized the extremal graphs which attain the minimum value for the distance Laplacian energy among all connected bipartite graphs and among all connected graphs with given independence number. The following problems will be of interest for the future research.\\
{\bf Problem 2.} Characterize the extremal graphs which attain the maximum value for the distance Laplacian energy among all connected bipartite graphs of order $n$.\\
{\bf Problem 3.} Characterize the extremal graphs which attain the maximum value for the distance Laplacian energy among all connected graphs of order $n$  with independence number $\alpha$.\\

\section{Distance Laplacian energy of graphs with given connectivity}

In this section, we obtain the distance Laplacian spectrum of the join of a connected graph $G_0$ with the union of two connected graphs $G_1$ and $G_2$. We show the existence of some new families of graphs having all the  distance Laplacian eigenvalues  integers. We also determine the graph with the minimum distance Laplacian energy among all the connected graphs with given vertex connectivity.

The vertex connectivity of a graph $ G $, denoted by $\kappa(G)  $, is the minimum number of vertices of $ G $ whose deletion disconnects $ G $. Let $ \mathcal{F}_n $ be the family of simple connected graphs on $ n $ vertices and let $ \mathcal{V}_n^k =\{G\in\mathcal{F}_n:\kappa(G)\leq k\}$.

Let $G_{1}(V_1,E_1)$ and $G_{2}(V_2,E_2)$ be two graphs on disjoint vertex sets $V_1$ and $V_2$ of order $n_1$ and $n_2$,  respectively. The {\em union} of graphs $G_1$ and $G_2$ is the graph $G_1\cup G_2=(V_1\cup V_2,E_1\cup E_2).$  The \textit{join} of graphs $G_1$ and $G_2$, denoted by $G_1\bigtriangledown G_2$ is the graph consisting of $G_1\cup G_2$ and all edges joining the vertices in  $V_1$ and the vertices in $V_2.$ In other words, the join of two graphs $G_{1}$ and $G_{2},$ denoted by $G_{1}\bigtriangledown G_{2},$ is the graph obtained from $G_{1}$ and $ G_{2}$ by joining each vertex of $G_{1}$ to every vertex of $G_{2}.$

In the following theorem, we find the distance Laplacian spectrum of the join of a connected graph $G_0$ with the union of two connected graphs $G_1$ and $G_2$ in terms of the distance Laplacian spectrum of the graphs $G_0,G_1$ and $G_2$.
\begin{theorem}\label{thm24}
Let $ G_0, G_1$ and $ G_2 $ be the connected graphs of order $ n_0, n_1  $ and $ n_2 $ respectively.  The distance Laplacian spectrum of $ G_0\bigtriangledown(G_1\cup G_2) $ consists of eigenvalues $ (\lambda+n_1+n_2), (\mu+n_0+2n_2), (\zeta+n_0+2n_1) ,2n-n_0,n,0$, where $ n=n_{0}+n_{1}+n_{2} $ and $ \lambda,\mu,\zeta $ varies over  the non-zero distance Laplacian eigenvalues of $ G_0, G_1 $, $ G_2  $, respectively.
\end{theorem}
\noindent\textbf{Proof.} Let $G_i$, $i=0,1,2$ be the connected graphs of order $n_i$. Let $ G= G_0\bigtriangledown(G_1\cup G_2) $ be the join of graphs $ G_0 $ and $ G_1\cup G_2 $. Clearly, $ G $ is a graph of diameter $ 2 .$ Let $D^L(G_0), D^L(G_1)$ and $D^L(G_2)$ be respectively the distance Laplacian matrices of the graphs $G_0,G_1$ and $G_2$.  By suitably labelling the vertices of $ G ,$ it can be seen that the distance Laplacian matrix of $ G $ is
\[D^{L}(G)=\begin{pmatrix}
b_0&-J_{n_0\times n_1}& -J_{n_0\times n_2}\\
-J_{n_1\times n_0}& b_{1}&-2J_{n_1\times n_2}\\
-J_{n_2\times n_0}&-2J_{n_2\times n_1}& b_{2}
\end{pmatrix}
,\] where $ b_0=(2n-n_{1}-n_{2})I_{n_{0}}-2J_{n_0}-D^L(G_0), b_{i}=(2n-n_{0})I_{n_{i}}-2J_{ n_i}-D^L(G_i) $, for $ i=1,2 $, $J_{n_{i}\times n_{j}}$ is the all one matrix of order $n_i\times n_j$, $I_{n_k}$ is the identity  matrix of order $n_k$. Let $ e_{n_i}=(1,1,\dots,1) $ be the all $ 1-$vector of order $ n_i, i=0,1, 2.$ It is well known that $ e_{n_i} $ is an eigenvector of $ G_i $ for the distance Laplacian  eigenvalue $ 0 $ and any other eigenvector $ x $ of $ G_i $ is orthogonal to $ e_{n_i} $. Let $ x $ be any eigenvector of $ G_0 $ for the non-zero distance Laplacian eigenvalue $ \lambda  $. Then $ x\perp e_{n_0}.$
Consider the column vector $ X=(x^{T}\ \ \ \ 0^{T}_{n_1\times 1}\ \ \ \ 0^{T}_{n_2\times 1})^T $, We have
\[
D^L(G)X =\begin{pmatrix}
(2n-n_{1}-n_{2})I_{n_{0}}-2J_{n_0}-D^L(G_0)x\\0\\0
\end{pmatrix}=(2n-n_{1}-n_{2}-\lambda) X
\]
giving that $ 2n-n_{1}-n_{2}-\lambda $ is an eigenvalue of $ D^L(G) $ for each non-zero distance Laplacian  eigenvalue $ \lambda $ of $ G_0 .$ In this way, we get $ n_0-1 $ distance Laplacian eigenvalues of $ G .$ Similarly, if $ 0\neq \mu $ is a distance Laplacian  eigenvalue of $ G_1 $ with eigenvector $ y$, $y \perp  e_{n_1}$, then it can be seen that the column vector $ Y=(0^{T}_{n_0\times 1}\ \ \ y^{T} \ \ \ 0^{T}_{n_2\times 1} )^T $ is an eigenvector of $ D^L(G) $ for the eigenvalue $ 2n-n_{0}-\mu$. This implies that $ 2n-n_{0}-\mu$ is an eigenvalue of $ D^L(G) $ for each non-zero distance Laplacian  eigenvalue $ \mu $ of $ G_1 .$ In this way, we get other $ n_1-1 $ distance Laplacian  eigenvalues of $ G.$ Lastly, if $ 0\neq \zeta $ is a distance Laplacian  eigenvalue of $ G_2 $ with eigenvector $ z, z\perp e_{n_2}$, then it can be seen that the column vector $ Z=(0^{T}_{n_0\times 1}\ \ \ 0^{T}_{n_1\times 1} \ \ \ z^{T} )^T $ is an eigenvector of $ D^L(G) $ for the eigenvalue $ 2n-n_{0}-\zeta$. This implies that $ 2n-n_{0}-\zeta$ is an eigenvalue of $ D^L(G) $ for each non-zero distance Laplacian  eigenvalue $ \zeta$ of $ G_2 .$ In this way, we get another $ n_2-1 $ distance Laplacian eigenvalues of $ G.$  Thus, we get a total of $ n_0-1+n_1-1+n_2-1=n-3 $ distance Laplacian  eigenvalues of $ G $. The remaining three distance Laplacian eigenvalues of $G$ are given by the quotient matrix of $D^{L}(G)$, which is
\[M=\begin{pmatrix}
n_1+n_2&-n_1&-n_2\\
-n_0&n_0+2n_2&-2n_2\\
-n_0&-2n_1&n_0+2n_1
\end{pmatrix}
.\] Since each row sum of  matrix $M$ is zero, it follows that one of the eigenvalues of this matrix is $0$. The other two eigenvalues of $M$ are given by the roots of $x^2-(3n-n_0)x+n(2n-n_0)=0$. This completes the proof.      \qed

A graph $G$ is said to an adjacency integral(Laplacian integral, distance integral) graph if all its adjacency eigenvalues(Laplacian eigenvalues, distance eigenvalues) are integers. Likewise, a graph $G$ is called distance Laplacian integral graph if all its distance Laplacian eigenvalues are integers.  From Theorem \ref{thm24} , it is clear that if all the graphs  $G_0$, $G_1$ and $G_2$ are distance Laplacian integral, then the graph $G= G_0\bigtriangledown(G_1\cup G_2) $ is also a distance Laplacian integral graph.  It is well known that the  complete graph $K_n$, the complete bipartite graph $K_{a,n-a}$, the complete split graph $CS_{a,n-a}$, the pineapple graph $PA_{n,p}$ (the graph obtained from a clique $K_{n-p}$ by adding $p$ pendent edges to a vertex of $K_{n-p}$) and the graph $S^{+}$(the graph obtained from the star graph $K_{1,n-1}$ by adding an edge between two pendent vertices) are all Laplacian integral graphs. Therefore, using Theorem \ref{thm24}, we have the following observation.
\begin{theorem}\label{thm25}
Let $n=n_0+n_1+n_2$, Then each of the graphs $K_{n_0}\bigtriangledown(K_{n_1}\cup K_{n_2})$, $K_{n_0}\bigtriangledown(K_{n_{1}-a,a}\cup K_{n_2})$, $K_{n_{0}-a,a}\bigtriangledown(K_{n_1}\cup K_{n_2})$, $K_{n_{0}-a,a}\bigtriangledown(K_{n_{1}-b,b}\cup K_{n_2})$, $K_{n_{0}-a,a}\bigtriangledown(K_{n_{1}-b,b}\cup K_{n_{2}-c,c})$, $CS_{a,n_{0}-a}\bigtriangledown(K_{n_1}\cup K_{n_2})$, $CS_{a,n_{0}-a}\bigtriangledown(CS_{b,n_{1}-b}\cup K_{n_2})$, $CS_{a,n_{0}-a}\bigtriangledown(CS_{b,n_{1}-b}\cup CS_{c,n_{2}-c}$, $K_{n_{0}-a,a}\bigtriangledown(CS_{b,n_{1}-b}\cup K_{n_2})$, $PA_{n_0,p}\bigtriangledown(K_{n_1}\cup K_{n_2})$, $CS_{a,n_{0}-a}\bigtriangledown(PA_{n_1,p}\cup K_{n_2})$, $CS_{a,n_{0}-a}\bigtriangledown(PA_{n_1,p}\cup K_{n_{2}-c,c})$,   etc are distance Laplacian integral graphs.
\end{theorem}

The following theorem shows that among all connected graphs with given vertex connectivity $k$ the graph $K_{k}\bigtriangledown(K_{t}\cup K_{n-k-t})$ has the minimum distance Laplacian energy.
\begin{theorem}\label{thm26}
Let $ G\in \mathcal{V}_n^k$ be a connected graph of order $ n\ge 4$ having vertex connectivity number $ k $. Then
$DLE(G)\ge t(2n-k-t+1)-\frac{2tW(G)}{n}$, for $1\le k<\frac{n-2t}{2}-\frac{n}{2t}$; $DLE(G)\ge t(3n-2k-2t)+n(n-t-k-1)-\frac{2(n-k-1)W(G)}{n}$, for $\frac{n-2t}{2}-\frac{n}{2t}\le k< n-t-\frac{n}{2t}$ and $DLE(G)\ge n(n-1)+2t(n-k-t)-\frac{2(n-1)W(G)}{n}$, for $n-t-\frac{n}{2t}\le k\leq n-2$; if $k=n-1$, then $G\cong K_n$ and so $DLE(G)=2n-2$. Equality occurs in each of these inequalities if and only if $G\cong K_{k}\bigtriangledown(K_{t}\cup K_{n-k-t}), 1\leq t \leq \lfloor\frac{n-k}{2}\rfloor $.
\end{theorem}
\noindent\textbf{Proof.}
Let $ G $ be a connected graph of order $ n $ with vertex connectivity $ k $, $1\leq k\leq n-1$. We first show
that $U^{L}_i(G)\geq U^{L}_i(K_{k}\bigtriangledown(K_{t}\cup K_{n-t-k}))$, for all $i=1,2,\dots,n$. If $k=n-1$, then $G\cong K_n$ and $K_{k}\bigtriangledown(K_{t}\cup K_{n-t-k})=K_n$ and so equality holds in this case. Assume that $1\leq k\leq n-2$, that is, $G$ is not a complete graph. Suppose that $ G$ is the connected graph of order $n$ with vertex connectivity $k$ for which the spectral parameter $U^{L}_{i}(G)$ has the minimum possible value. Then, it is clear that $ G\in \mathcal{V}_{n}^{k}$ and $U^{L}_{i}(G)$ attains the minimum value for $G$. Let $S$ be a vertex cut set of $G$ with $|S|=k$. Let $G_1,G_2,\dots,G_r$ be the connected components of the graph $G-S$. We will show the number of components of graph $G-S$ is two, that is $r=2$. For if, $r>2$, then adding an edge between any two components, say $G_1$ and $G_2$ of $G-S$ gives the graph $G^{'}=G+e$, which is such that the vertex connectivity of $G^{'}$ is $k$. Clearly, $G^{'}\in\mathcal{V}_n^k$, also by Lemma \ref{lem1}, we have $ U^{L}_i(G)> U^{L}_i(G^{'})$. This is a contradiction to the fact $U^{L}_i(G)$ attains the  minimum possible value  for $G$. Therefore, we must have $r=2$. Further, we claim that each of  the components $G_1$, $G_2$ and the vertex induced subgraph $\langle S\rangle$ are cliques. For if one among them say $G_1$ is not a clique, then adding an edge between the two non adjacent vertices of $G_1$ gives a graph $H$ having vertex connectivity same as the vertex connectivity of $G$. Clearly, $H\in\mathcal{V}_n^k$ and by Lemma \ref{lem1}, we have $ U^{L}_i(G)> U^{L}_i(H)$, Which is a contradiction as $U^{L}_i(G)$ attains the  minimum possible value  for $G$. Thus, $G$ must be of the form $G= K_{k}\bigtriangledown(K_{t}\cup K_{n-k-t}), 1\leq t \leq \lfloor\frac{n-k}{2}\rfloor $. This shows that for all $ G\in \mathcal{V}_n^k$, the  spectral parameter $U^{L}_i(G)$ has the minimum possible value for the graph $K_{k}\bigtriangledown(K_{t}\cup K_{n-k-t})$. That is, for all $ G\in \mathcal{V}_n^k$, we have
$U^{L}_i(G)\geq U^{L}_i(K_{k}\bigtriangledown(K_{t}\cup K_{n-t-k})).$ With this, from the definition of distance Laplacian energy, it follows that
 \begin{align}\label{33}
 DLE(G)&=2\left(U_{\sigma}^L(G) -\frac{2\sigma W(G)}{n}\right) =2\max_{1\leq j\leq n}\left( \sum\limits_{i=1}^{j}\rho_i^L(G)-\frac{2jW(G)}{n} \right)\nonumber\\&
 \geq 2\max_{1\leq j\leq n}\left( \sum\limits_{i=1}^{j}\rho_i^L(K_{k}\bigtriangledown(K_{t}\cup K_{n-t-k}))-\frac{2jW(G)}{n} \right).
 \end{align}
Taking $G_{n_0}=K_k,G_{n_1}=K_t, G_{n_2}=K_{n-k-t}, n_0=k,n_1=t$ and $n_2=n-t-k$ in Theorem \ref{thm24}, we find that the distance Laplacian spectrum of the graph $K_{k}\bigtriangledown(K_{t}\cup K_{n-t-k})$ is $\{2n-k,(2n-t-k)^{[t-1]},(n+t)^{[n-t-k-1]},n^{[k]},0\}$. Let $ \sigma$ be the number of distance Laplacian eigenvalues of $K_{k}\bigtriangledown(K_{t}\cup K_{n-t-k})$ which are greater than or equal to that $\frac{2W(K_{k}\bigtriangledown(K_{t}\cup K_{n-t-k}))}{n}=\frac{n^2-n+2nt-2t^2-2kt}{n}.$ It is easy to see that $2n-k$ is the distance Laplacian spectral radius of the graph $K_{k}\bigtriangledown(K_{t}\cup K_{n-t-k})$ and so for this eigenvalue, we always have $2n-k\geq \frac{2W(K_{k}\bigtriangledown(K_{t}\cup K_{n-t-k}))}{n}$. For the eigenvalue $2n-k-t$, we have $2n-k-t\geq \frac{2W(K_{k}\bigtriangledown(K_{t}\cup K_{n-t-k}))}{n}=\frac{n^2-n+2nt-2t^2-2kt}{n}$ giving that
\begin{align}\label{d}
2t^2-(3n-2k)t+(n^2+n-kn)\ge 0.
\end{align}
The zeros of the polynomial $g_1(t)=2t^2-(3n-2k)t+(n^2+n-kn)$ are $y_1=\frac{3n-2k+\sqrt{(n-2k)^2-8n}}{2}$
 and $y_2=\frac{3n-2k-\sqrt{(n-2k)^2-8n}}{2}$. This shows that $g_1(t)\ge 0$, for all $t\le y_2$ and $t\ge y_1$. Since, $\frac{n-k}{2}< \frac{3n-2k-\sqrt{(n-2k)^2-8n}}{2}=y_2$ always holds, it follows that $g_1(t)\ge 0$, for all $t\le \frac{n-k}{2}$. For the eigenvalue $n+t$, we have $n+t\geq \frac{2W(K_{k}\bigtriangledown(K_{t}\cup K_{n-t-k}))}{n}=\frac{n^2-n+2nt-2t^2-2kt}{n}$ giving that $2t^2-(n-2k)t+n\ge 0$, which in turn gives $k\geq \frac{n-2t}{2}-\frac{n}{2t}$. This shows that $n+t\geq \frac{2W(K_{k}\bigtriangledown(K_{t}\cup K_{n-t-k}))}{n}$, for all $k\geq \frac{n-2t}{2}-\frac{n}{2t}$ and $n+t< \frac{2W(K_{k}\bigtriangledown(K_{t}\cup K_{n-t-k}))}{n}$, for all $k< \frac{n-2t}{2}-\frac{n}{2t}$.  Lastly, for the eigenvalue $n$, we have $n\ge  \frac{2W(K_{k}\bigtriangledown(K_{t}\cup K_{n-t-k}))}{n}=\frac{n^2-n+2nt-2t^2-2kt}{n}$ giving that
   $2t^2-(2n-2k)t+n\ge  0$, which in turn gives that $k\geq n-t-\frac{n}{2t}$. This shows that $n\geq \frac{2W(K_{k}\bigtriangledown(K_{t}\cup K_{n-t-k}))}{n}$, for all $k\geq n-t-\frac{n}{2t}$ and $n< \frac{2W(K_{k}\bigtriangledown(K_{t}\cup K_{n-t-k}))}{n}$, for all $k< n-t-\frac{n}{2t}$. From this discussion it follows that, if $1\leq k< \frac{n-2t}{2}-\frac{n}{2t}$, then $\sigma=t$, if $\frac{n-2t}{2}-\frac{n}{2t}\le k<n-t-\frac{n}{2t}$, then $\sigma=n-k-1$ and if $k\ge n-t-\frac{n}{2t}$, then $\sigma=n-1$. If $1\le k<\frac{n-2t}{2}-\frac{n}{2t}$, then from (\ref{33}), it follows that
\begin{align*}
 DLE(G)&\geq 2\max_{1\leq j\leq n}\left( \sum\limits_{i=1}^{j}\rho_i^L(K_{k}\bigtriangledown(K_{t}\cup K_{n-t-k}))-\frac{2jW(G)}{n} \right)\\&
 \geq 2\left( \sum\limits_{i=1}^{t}\rho_i^L(K_{k}\bigtriangledown(K_{t}\cup K_{n-t-k}))-\frac{2tW(G)}{n} \right)\\&
 =t(2n-k-t+1)-\frac{2tW(G)}{n}.
\end{align*}
If $\frac{n-2t}{2}-\frac{n}{2t}\le k< n-t-\frac{n}{2t} $, then from (\ref{33}), it follows that
\begin{align*}
 DLE(G)&\geq 2\max_{1\leq j\leq n}\left( \sum\limits_{i=1}^{j}\rho_i^L(K_{k}\bigtriangledown(K_{t}\cup K_{n-t-k}))-\frac{2jW(G)}{n} \right)\\&
 \geq 2\left( \sum\limits_{i=1}^{n-k-1}\rho_i^L(K_{k}\bigtriangledown(K_{t}\cup K_{n-t-k}))-\frac{2(n-k-1)W(G)}{n} \right)\\&
 =t(3n-2k-2t)+n(n-t-k-1)-\frac{2(n-k-1)W(G)}{n}.
\end{align*}
If $n-t-\frac{n}{2t}\le k\le n-1$, then from (\ref{33}), it follows that
\begin{align*}
 DLE(G)&\geq 2\max_{1\leq j\leq n}\left( \sum\limits_{i=1}^{j}\rho_i^L(K_{k}\bigtriangledown(K_{t}\cup K_{n-t-k}))-\frac{2jW(G)}{n} \right)\\&
 \geq 2\left( \sum\limits_{i=1}^{n-1}\rho_i^L(K_{k}\bigtriangledown(K_{t}\cup K_{n-t-k}))-\frac{2(n-1)W(G)}{n} \right)\\&
 =n(n-1)+2t(n-k-t)-\frac{2(n-1)W(G)}{n}.
\end{align*}
This completes the proof.\qed

The next observation  follows from the proof of the Theorem \ref{thm26}.
\begin{corollary}
Let $G$ be a connected graph of order $n\ge 4$ having  vertex connectivity $\kappa\le k$. Let $\rho_{1}^{L}(G)\ge \rho_{2}^{L}(G)\ge\cdots\ge \rho_{n-1}^{L}(G)>\rho_{n}^{L}(G)=0$ be the distance Laplacian eigenvalues of $G$. Then $\rho_{1}^{L}(G)\geq 2n-k,$ $\rho_{i}^{L}(G)\ge 2n-t-k$, for $2\le i\le t$, $\rho_{i}^{L}(G)\ge n+t$, for $t+1\le i\le n-k-1$ and $\rho_{i}^{L}(G)\ge n$, for $n-k\le i\le n-1$.
Equality occurs in each of these inequalities if and only if $G\cong K_{k}\bigtriangledown(K_{t}\cup K_{n-t-k}).$
\end{corollary}

In the later part of this section, we characterized the extremal graphs which attain the minimum value for the distance Laplacian energy among all connected graphs with given vertex connectivity. The following problem will be of interest for the future research.\\
{\bf Problem 4.} Characterize the extremal graphs which attain the maximum value for the distance Laplacian energy among all connected graphs of order $n$ with given vertex connectivity.

\section{Conclusion}
Let $\mathbb{M}_n(\mathbb{C})$ be the set of all square matrices of order $n$ with complex entries. The trace norm of a matrix $M\in \mathbb{M}_n(\mathbb{C})$ is defined as $\lVert M\rVert_{*}=\sum\limits_{i=1}^{n}\sigma_i(M)$, where $\sigma_1(M)\ge \sigma_2(M)\ge\cdots\ge \sigma_n(M)$ are the singular values of $M$(i.e. the square roots of the eigenvalues of $MM^*$, where $M^*$ is the complex conjugate of $M$). It is well known that for a symmetric matrix $M$, if $\sigma_i(M)$ the singular values and  $\lambda_i(M)$, $ i=1,2,\dots,n,$  are the eigenvalues, then $ \sigma_i(M)=|\lambda_i(M)|$. In the light of this definition, it follows that the distance Laplacian energy $DLE(G)$ of a connected graph $G$ is the trace norm of the matrix $ D^L(G)-\frac{2W(G)}{n}I_n$, where $I_n$ is the identity matrix of order $n$. It is an interesting problem in Matrix theory to determine among a given  class of matrices the matrix (or the matrices) which attain the maximum value and the minimum value for the trace norm. The trace norm of matrices associated with the graphs and digraphs are extensively studied. For some recent papers in this direction see  \cite{mr} and the references therein. \\
\indent Therefore, in this language Theorem \ref{thm23} gives relation between the trace norm of the matrices $ D^L(G)-\frac{2W(G)}{n}I_n$ and  $ L(G)-\frac{2m}{n}I_n$, when $G$ is a connected graph of diameter two;  Theorem \ref{thm2} gives a lower bound for the trace norm of $ D^L(G)-\frac{2W(G)}{n}I_n$ in terms of the order and the trace of the matrix $D^L(G)$; Theorem \ref{thm31} gives that among all bipartite connected graphs $G$, the complete bipartite $K_{a,b}$ attains the minimum trace norm for the matrix $D^L(G)-\frac{2W(G)}{n}I_n$; Theorem \ref{thm3} gives that among all connected graphs $G$ with given independence number $n-t$, $1\le t\le n-1$, the complete split graph $CS_{t,n-t}$ attains the minimum trace norm for the matrix $D^L(G)-\frac{2W(G)}{n}I_n$ and Theorem \ref{thm26} gives that among all connected graphs $G$ with given vertex connectivity $k$ the graph $K_{k}\bigtriangledown(K_{t}\cup K_{n-k-t})$ attains the minimum trace norm for the matrix $D^L(G)-\frac{2W(G)}{n}I_n$.

\noindent{\bf Acknowledgements.} The research of S. Pirzada is supported by SERB-DST, New Delhi under the research project number MTR/2017/000084.


\begin{thebibliography}{0}


\bibitem{AH1} {M. Aouchiche and  P. Hansen,  Distance spectra of graphs: a survey,  \em Linear Algebra Appl.} {\bf 458} (2014) 301--386.
\bibitem{AH2} {M. Aouchiche and  P. Hansen,  Two Laplacians for the distance matrix of a graph, \em Linear Algebra Appl.} {\bf 439} (2013) 21--33.
\bibitem{AH} {M. Aouchiche, P. Hansen,  Some properties of the distance Laplacian eigenvalues of a graph,  \em Czechoslovak Mathematical journal}, {\textbf{64}} (139)  (2014) 751-761.

\bibitem{DMCM} {D. M. Cvetkovi\'{c}, M. Doob and H. Sachs, \em Spectra of graphs. Theory and application,} Pure and Applied Mathematics, 87. Academic Press, Inc. New York, 1980.

\bibitem{dah} {K. C. Das,  M. Aouchiche and P. Hansen, On (distance) Laplacian energy and (distance) signless Laplacian energy of graphs, \em Discrete Appl. Math.} {\bf 243} (2018) 172--185.

\bibitem{dr} {R. C. Diaz and O. Rojo, Sharp upper bounds on the distance energies of a graph, \em Linear Algebra Appl.} {\bf 545} (2018) 55--75.

\bibitem{hilal} {H. A. Ganie, On the distance Laplacian spectrum(energy) of graphs, \em Discrete Math. Algorithms Appl.}, 2050061 (2020), (16 pages).

\bibitem{hbp} {H. A. Ganie, B. A. Chat and S. Pirzada, On the signless Laplacian energy of a graph and energy of line graph, \em Linear Algebra Appl.} {\bf 544} (2018) 306--324.

\bibitem{hpbt} {Hilal A. Ganie, S. Pirzada, B. A. Rather and V. Trevisan, Further developments on Brouwer's conjecture for the sum of Laplacian eigenvalues of graphs, \em Linear Algebra Appl.}, {\bf 588} (2020), 1-18.

\bibitem{IGBZ} {I. Gutman and B. Zhou,  Laplacian energy of a graph, \em Linear Algebra Appl.} {\bf 414} (2006) 29--37.

\bibitem{IGV} {G. Indulal, I. Gutman and A. Vijayakumar,  On distance energy of graphs, \em MATCH Commun. Math. Comput. Chem.} {\bf 60} (2008) 461--472.

\bibitem{lsg} {X. Li, Y. Shi and I. Gutman, \em Graph Energy}, Springer, New York, 2012.
\bibitem{mr} {J. Monsalve and  J. Rada, Oriented bipartite graphs with minimal trace norm, \em Linear Multilinear Algebra}, {\bf 67} (2019), 1121-1131.

\bibitem{ph} {S. Pirzada and H. A. Ganie, On the Laplacian eigenvalues of a graph and Laplacian energy, \em Linear Algebra Appl.} {\bf 486} (2015) 454--468.
\bibitem{sp} {S. Pirzada, \em An Introduction to Graph Theory}, Universities Press, Orient BlackSwan, Hyderabad (2012).
\bibitem{YYG} {J. Yang, L. You and I. Gutman,  Bounds on the distance Laplacian energy of graphs, \em Kragujevac J. Math.} {\bf 37} (2013) 245--255.


\end{thebibliography}
\end{document}